\documentclass[12pt]{amsart}

\usepackage[hmargin=2cm,bmargin=2.5cm,tmargin=3cm]{geometry}

\usepackage[utf8]{inputenc}
\usepackage{eucal,mathrsfs}

\usepackage[usenames]{color}

\usepackage[normalem]{ulem}
\usepackage{epigraph}	

\usepackage{amsmath,amssymb,amsthm,amsfonts,enumerate,url,mathbbol,mathrsfs,tikz,graphicx,pifont,doi,comment}

\usepackage{tikz, subfigure, xcolor} 
\usepackage{pgfplots}
\usepackage{cite}

\usepackage{dsfont}

\usepackage{hyperref}

\newtheorem{theorem}{Theorem}[section]

\theoremstyle{definition}
\newtheorem{defi}[theorem]{Definition}
\newtheorem{remark}[theorem]{Remark}
\newtheorem{example}[theorem]{Example}

\usepackage{faktor}

\numberwithin{equation}{section}

\DeclareMathOperator{\Graph}{\mathcal G}
\DeclareMathOperator{\Graphprim}{\mathfrak G}

\newcommand{\R}{\mathbb{R}}
\newcommand{\K}{\mathbb{K}}

\newcommand\restr[2]{{% we make the whole thing an ordinary symbol
 \left.\kern-\nulldelimiterspace % automatically resize the bar with \right
 #1 % the function
 \vphantom{\big|} % pretend it's a little taller at normal size
 \right|_{#2} % this is the delimiter
 }}

\usepackage{eucal}

 \def\mG{\mathsf{G}}
 \def\mV{\mathsf{V}}
 \def\mE{\mathsf{E}}

 \def\mv{\mathsf{v}}
 \def\me{\mathsf{e}}
 \def\mw{\mathsf{w}}
 \def\mf{\mathsf{f}}

\author{Delio Mugnolo}
\address{Chair of Analysis, Faculty of Mathematics and Computer Science, University of Hagen, 58084 Hagen, Germany}
\email{delio.mugnolo@fernuni-hagen.de}
\date{\today}
\thanks{I would like to thank Amru Hussein, James Kennedy, Pavel Kurasov, and Marvin Plümer for useful comments.}

\title[What is actually a metric graph?]{What is actually a metric graph?}

\begin{document}

\begin{abstract}
Metric graphs are often introduced based on combinatorics, upon ``associating'' each edge of  a graph with an interval; or else, casually ``gluing'' a collection of intervals at their endpoints in a network-like fashion. Here we propose an abstract, self-contained definition of metric graph. Being mostly topological, it doesn't require any knowledge from graph theory and already determines uniquely several concepts that are commonly and unnecessarily \textit{defined} in the literature. Nevertheless, many ideas mentioned here are folklore in the quantum graph community: we discuss them for later reference. 
%They can easily be extended to related settings, like hypergraphs and simplicial complexes.
\end{abstract}

\maketitle

\section{Metric graphs as quotient spaces}
Let $\mE$ be a countable set. Given some $(\ell_\me)_{\me\in\mE}\subset (0,\infty)$, we consider the family $[0,\ell_\me]_{\me\in\mE}$ of metric measure subspaces of $\R$ (wrt Euclidean metric $d_\me$ and Lebesgue measure $\lambda_\me$) and their disjoint union
\[
\mathcal E:=\bigsqcup\limits_{\me\in\mE} [0,\ell_\me]\ :
\]
we adopt the usual notation $(x,\me)$ for the element of $\mathcal E$ with $x\in [0,\ell_\me]$ and $\me\in\mE$.

\medskip
We endow $\mathcal E$ with the \textit{disjoint union topology}: by definition, this means that a subset $U$ of $\mathcal E$ is open if and only if its preimage $\varphi _{\me}^{-1}(U)$ is open in $[0,\ell_\me]$ for each $\me\in\mE$, where $\varphi_\me$ is the canonical injection $\varphi_\me:[0,\ell_\me]\ni x\mapsto (x,\me)\in \mathcal E$. Hence a set $U$ is open if and only if each $\varphi _{\me}^{-1}(U)$ is a union of sets of the form $[0,\varepsilon_1)$, $(\varepsilon_2,\ell_\me]$, or $(\varepsilon_3,\varepsilon_4)$, for $\varepsilon_i\in (0,\ell_\me)$.  Disjoint unions of such sets thus form a basis of the topology of $\mathcal E$.

 The disjoint union topology  of $\mathcal E$ is metrizable and indeed it agrees with the topology induced by the metric defined by setting
\begin{equation}\label{eq:de}
d_{\mathcal E}\big((x,\me),(y,\mf)\big):=\begin{cases}
d_\me(x,y)=|x-y|,\qquad &\hbox{if }\me=\mf\hbox{ and }x,y\in [0,\ell_\me],\\
\infty, &\hbox{otherwise}.
\end{cases}
\end{equation}
(We follow throughout the convention in~\cite{BurBurIva01} and adopt the generalized notion of distance that allows for the value $\infty$.)

Consider the set 
\[
\mathcal V:=\bigsqcup\limits_{\me\in\mE}\{0,\ell_\me\}
\]
of endpoints of $\mathcal E$. 
Given any equivalence relation $\sim$ on $\mathcal V$, we extend it to an equivalence relation on $\mathcal E$ by equality: i.e., two elements $(x_1,\me_1),(x_2,\me_2)\in \mathcal E$ belong to the same equivalence class if and only if $(x_1,\me_1)=(x_2,\me_2)$ or else $(x_1,\me_1),(x_2,\me_2)\in \mathcal V$ and $(x_1,\me_1)\sim(x_2,\me_2)$. With an abuse of notation we denote this equivalence relation on $\mathcal E$ again by $\sim$: this allows us to introduce quotient sets.

\begin{defi}\label{def:metricg}
We call $\Graph:=\faktor{\mathcal E}{\sim}$ a \textit{metric graph} and $\mV:=\faktor{\mathcal V}{\sim}$ its \textit{set of vertices}. 
\end{defi}

 Let us stress that $\mathcal G$ is, canonically, a topological space with respect to the quotient topology, see Remark~\ref{rem:quot}. We will devote part of Section~\ref{sec:functsp} to the study of functions that are continuous over this topological space.
 
\begin{remark}
This setting can be slightly generalized by considering an additional countable set $\mE_\infty$ and replacing the sets $\mathcal E$ and $\mathcal V$ studied above by
\[
\mathcal E:=\bigsqcup_{\me\in \mE}[0,\ell_\me]\sqcup \bigsqcup_{\me\in \mE_\infty}[0,\infty)
\]
and
\[
\mathcal V:=\bigsqcup_{\me\in \mE}\{0,\ell_\me\}\sqcup \bigsqcup_{\me\in \mE_\infty}\{0\},
\]
respectively. In this way we can add semi-infinite leads to a metric graph consisting of a ``countable core of bounded edges''.
\end{remark}

According to this definition, a metric graph is uniquely determined by a family $(\ell_\me)_{\me\in \mE}$ and an equivalence relation on $\mathcal V$. Its vertices are the cells of the partition of   $\mathcal V$ induced by $\sim$. Two vertices $\mv,\mw\in \mV$ are said to be \textit{adjacent} if there exists some (not necessarily unique) $\me\in\mE$ such that $\{x,y\}=\{0,\ell_\me\}$ for  representatives $x$ of $\mv$ and $y$ of $\mw$; in this case we write $\mv\sim \mw$ and, with an abuse of notation, also $\mv\sim \me$. The cardinality $\deg(\mv)$ of the set $\{\mw\in\mV: 
\mw\hbox{ is adjacent to }\mv\}$ is called \textit{degree} of $\mv\in \mV$; $\Graph$ is called \textit{combinatorially locally finite} if $\deg(\mv)<\infty$ for all $\mv\in\mV$, and \textit{metrically locally finite} if $\sum\limits_{\substack{\me\in\mE\\ \mv\sim \me}} \ell_\me<\infty$ for all $\mv\in\mV$.

\medskip
To justify Definition~\ref{def:metricg}, we are going to show how $\Graph$ can be canonically endowed with a metric.
Following \cite[Def.~3.1.12]{BurBurIva01} we introduce the quotient pseudo-metric defined by
\[
d_{\Graph}(\xi,\theta):=\inf \sum_{i=1}^k d_{\mathcal E}(\xi_i,\theta_i),\qquad \xi,\theta\in \Graph,
\]
where the infimum is taken over all $k\in \mathbb N$ and all pairs of $k$-tuples $(\xi_1,\ldots, \xi_k)$ and $(\theta_1,\ldots, \theta_k)$ with $\xi=\xi_1$, $\theta=\theta_k$, and $\theta_i\sim \xi_{i+1}$ for all $i=1,\ldots,k-1$. We call $d_{\Graph}$ the \textit{path pseudo-metric} of $\Graph$.

\begin{defi}
A metric graph  is \textit{connected} if the path pseudo-metric  doesn't attain the value $\infty$. 
\end{defi}

\begin{remark}
While $d_{\Graph}$ is \textit{a priori} only a pseudo-metric, it is actually a (generalized) metric (i.e., $d_{\Graph}(\xi,\theta)=0$ implies $\xi=\theta$; but the value $\infty$ can still be attained), which we call the \textit{path metric} of $\Graph$, if $\mE$ is finite or, more generally, if $\inf\limits_{\me\in\mE}\ell_\me>0$. 

Alternatively, consider the \textit{doubly connected part} $\mathcal G_d$ of $\mathcal G$, i.e., the set of all $(x,\me)$, $x\in (0,\ell_\me)$, whose removal doesn't turn $\Graph$ into a disconnected metric graph. Let us assume $\mathcal G_d\ne \emptyset$ and denote by $\mE_d$ the set of its edges. Then
 $d_{\Graph}$ is a metric if $\inf\limits_{\me\in\mE_d}\ell_\me>0$.
\end{remark}

A connected metric graph is hence a metric space. Furthermore,
$d_{\mathcal E}$ is the \textit{disjoint union length (pseudo-)metric}, in the sense of \cite[Def.~3.1.15]{BurBurIva01}; hence $d_{\Graph}$ is actually a \textit{length (pseudo-)metric}, thus any connected metric graph is a \textit{length metric space} (in the sense of~\cite{Stu06}); and even a \textit{geodesic space} (again in the sense of~\cite{Stu06}) whenever $\inf\limits_{\me\in\mE}\ell_\me >0$. In the latter case, the metric space $\Graph$ is also complete, hence a Polish space.

Also, observe that because each interval is an Alexandrov space of nonnegative (actually, vanishing) curvature, any (possibly infinite) metric star is a bouquet (in the sense of~\cite[Def.~4.2.7]{BurBurIva01}) thereof, hence an Alexandrov space of nonpositive curvature by~\cite[Prop.~4.2.9]{BurBurIva01}: because each finite metric graph can be recursively obtained as a bouquet of metric stars, one concludes that each metric graph is an Alexandrov space of nonnegative curvature. In the language of~\cite{BurGroPer92}, the \textit{space of directions} can only be a singleton at vertices of degree 1: accordingly, the \textit{singular set} $S_{\Graph}$ is discrete.

\smallskip
The topology of $\mathcal G$ induced by the pseudo-metric $d_{\Graph}$ is easily described: 
%that the open balls wrt to the pseudometric $d_{\Graph}$ induce a topology on $\Graph$. It is immediate that 
a basis of this topology consists of open balls wrt to $d_{\Graph}$, i.e., of sets that are either open subintervals of $[0,\ell_\me]$'s (``open subsets of edges'') or -- up to gluing wrt $\sim$ -- disjoint unions of semi-open subintervals of $[0,\ell_\me]$'s (``open stars centred at vertices'').

\begin{remark}
We have already pointed out that $\Graph$ canonically becomes a topological space whenever endowed with the quotient topology\footnote{ I.e., a subset of $\Graph$ is open if and only if it consists of equivalence classes whose union is open in $\mathcal E$; it also follows that a subset of $\Graph$ is closed if and only if it consists of equivalence classes whose union is closed in $\mathcal E$.}. A basis of the topology of $\Graph$ is then given by images under the canonical surjection of elements of a basis of $\mathcal E$: in particular, disjoint unions of open subintervals of $[0,\ell_\me]$'s (``open subsets of edges'') and -- up to gluing wrt $\sim$ -- disjoint unions of semi-open subintervals of $[0,\ell_\me]$'s (``open stars centred at vertices'')\footnote{ It is known that given $q:\mathcal E\to \Graph$ and given a basis $\mathcal B$ of the topology of $\mathcal E$, $q(\mathcal B)$ is a basis of the topology of $\Graph$ if and only if $q$ is open, which is of course especially the case if $q$ is the canonical surjection.
%: \url{https://math.stackexchange.com/questions/2082488/basis-for-the-quotient-topology}.
}. Hence,  the canonical quotient topology on $\Graph$ coincides with the topology induced by the path pseudo-metric on $\Graph$.
\end{remark}

\medskip
Finally, $\Graph$ is clearly a measure space with respect to the direct sum measure $\mu=\bigoplus\limits_{\me\in\mE} \lambda_\me$~\cite[214K]{Fre03}; this measure space is finite if $\Graph$ has \textit{finite volume}, i.e., if $\mu(\Graph)=\sum\limits_{\me\in\mE} \ell_\me<\infty$. A sufficient condition for $\Graph $ to be  a metric measure space in the sense of~\cite[\S~3]{Stu06} is that $\inf\limits_{\me\in\mE}\ell_\me>0$.

\begin{remark}\label{rem:quot}
Let $\iota:=q\circ \partial$, where $\partial :\mE\ni \me\mapsto (0,\ell_\me)\in\mathcal V^2$ and $q$ is the canonical extension to $\mathcal V^2$ of the canonical surjection $q:\mathcal V\to \mV$  defined by $q(x,y):=\{q(x),q(y)\}$: the latter set may thus consist of either one or two elements of $\mV$. 

Then the triple $\mG:=(\mV,\mE,\iota)$ is a multigraph (recall that, by definition, a multigraph may have loops and parallel edges, s.~\cite[\S~1.10]{Die05}): we call it the  \textit{combinatorial multigraph underlying} $\Graph$.   The (pseudo)metric on $\Graph$ induces a (pseudo)metric on $\mG$ -- in fact, the canonical one commonly used in graph theory. We stress that, however, the edge lengths need not satisfy the triangle inequality: an example being given by a triangle with edge lengths 1, 2, and 4.
\end{remark}

\section{Graph surgery}

\begin{defi}
Let $\Graph$ be a metric graph. 
Let $\approx$ be a further  equivalence relation on $\mathcal V$. Then 
$\hat{\Graph}:=\faktor{\mathcal E}{\approx}$ is called a \textit{rewiring} of $\mathcal G=\faktor{\mathcal E}{ \sim}$; and a \textit{cut} (resp., \textit{non-trivial cut}) of $\mathcal G$ if $\approx$ is coarser (resp., strictly coarser) than $\sim$.
\end{defi}

Any function $f:\Graph\to \K$ canonically induces a function $\hat{f}:\mathcal E\to \K$, or equivalently  a family $\hat{f}=(f_\me)_{\me\in \mE}$ with $f_\me:[0,\ell_\me]\to \K$ for all $\me\in\mE$ (and hence on each rewiring of $\Graph$, and especially on each of its cuts); the converse is wrong, though, since the boundary values of $f$ in $\mathcal V$ may conflict with the equivalence relation $\sim$ that defines $\Graph$.
% More generally, any function $f:\Graph\to \K$ canonically induces a function on each of the cuts of $\Graph$, but the converse is generally wrong.

\medskip
We regard a cut $\hat{\Graph}$ as a new metric graph obtained by cutting through some vertices of $\Graph$. By definition, any cut of $\mathcal G$ shares with $\mathcal G$ its edge set $\mathcal E$: this can be limiting in certain situations and suggests to introduce the following.
%Among other things, cuts are interesting because they allow us to partition $\Graph$. The problem is that once $\mE$ is given, any cut of $\Graph$ defined in this way can only have at most $|\mE|$ connected components. For this reason, it is crucial that subdividing existing edges is allowed for.

\begin{defi}
Let $\Graph$ be a metric graph and $\hat{\mE}$ be a countable set  such that there exists a surjection $\varsigma:\hat{\mE}\to \mE$. Given a vector $(\ell_{\me})_{\me\in\hat{\mE}}$ and an equivalence relation $\hat{\sim}$ on $\hat{\mathcal V}=\bigsqcup\limits_{\me\in\hat{\mE}}\{0,\ell_{\me}\}$, consider the natural extension of $\hat{\sim}$ wrt  an induced surjection $\varsigma:\hat{\mathcal E}:=\bigsqcup\limits_{\me\in\hat{\mE}}[0,\ell_\me]\to \mathcal E$: given $x,y\in \hat{\mathcal E}$, $x\hat{\sim} y$ if  $\varsigma(x)=\varsigma(y)$.

Then the metric graph $\hat{\Graph}:=\faktor{\hat{\mathcal E}}{\hat{\sim}}$ is called a \textit{subdivision} of $\Graph$ if for all $\me\in\mE$ the set $\varsigma^{-1}(\me)$ can be enumerated in such a way, say $\varsigma^{-1}(\me)=\{\me_1,\ldots,\me_{k_\me}\}$, that
\begin{itemize}
\item $(0,\me)=(0,\me_1),\ (\ell_{\me_1},\me_1)\hat{\sim}(0,\ell_{\me_2}),\ldots, (\ell_{k_\me-1},{k_\me-1})\hat{\sim} (0,{k_\me}),\ (\ell_{k_\me},k_\me)=(\ell_\me,\me)$
\item $\sum\limits_{j=1}^{k_\me}\ell_{\me_j}=\ell_\me$.
\end{itemize}
\end{defi}

Given a connected metric graph $\Graph$, any two subdivisions of $\Graph$ are isometric metric spaces.

Roughly speaking, $\varsigma(\hat{\me})=\me$ if $\hat{\me}=\me_i$ for some $i\in\{1,\ldots,k_\me\}$ (i.e., if  $\me$ is the edge that has been split to produce $\me_1,\ldots,\me_{k_\me}$, one of which is precisely $\hat{\me}$); and $x\hat{\sim} y$ if $x,y$ are representatives of a new vertex that has been created in $\hat{\Graph}$ inside the edge $[0,\ell_\me]$. Again, each function $f=\bigoplus\limits_{\me\in\mE}f_\me:\Graph\to \K$ canonically induces a new function $\hat{f}=\bigoplus\limits_{\me\in\varsigma^{-1}(\mE)}f_{\me}$ on any subdivision $\hat{\Graph}$, but the converse is generally wrong.

\begin{remark}
Given a subdivision $\Graph'$ of a metric graph $\Graph$, the equivalence relation $\sim$ that defines $\Graph$ can be canonically identified with the the equivalence relation $\sim'$ that defines $\Graph'$; hence, the set of all equivalence relations on $\mathcal V$ can be canonically embedded in the class of all equivalence relations on $\mathcal V'$. Therefore, given two different subdivisions $\Graph',\Graph''$ of $\Graph$, there is always a new subdivision whose vertex set contains all vertices of both $\Graph'$ and $\Graph''$ (this defines a partial ordering on the set of subdivisions of $\Graph$).
\end{remark}

 In the literature,  surgery of metric graphs has been frequently performed according to these rules: metric graphs arising by cutting through vertices of $\Graph$ in the sense of~\cite[Def.~3.2]{BerKenKur19} are non-trivial cuts of subdivisions of $\Graph$, in the language of the present note; whereas metric graphs arising by transplantation (and especially unfolding) as in~\cite[Def.~3.15 and Def.~3.16]{BerKenKur19} are rewirings of subdivisions of $\Graph$. (Non-trivial) symmetrisations of edges  in~\cite[Def.~3.17]{BerKenKur19}, on the other hand, can not be described in terms of a subdivision's rewirings or cuts.

\begin{defi}
Let $\Graph$ be a metric graph. We call any metric graph  arising from a rewiring or cut of a subdivision of $\Graph$ as a  \textit{rearrangement} of $\Graph$. 
\end{defi}

While comparing rearrangements $\Graph_1,\Graph_2$ of $\Graph$ we can certainly assume without loss of generality that they are rewiring or cuts of the \textit{same} subdivision $\hat{\Graph}=\faktor{\hat{\mathcal E}}{\hat{\sim}}$ of $\Graph$.
%Let $\mE$ be a countable set and $(\ell_\me)_{\me\in\mE}\subset (0,\infty)$. 
While rearrangements of $\Graph$ generally have a different metric, given any two equivalence relations $\approx_1$ and $\approx_2$ on $\hat{\mathcal V}$ and the associated canonical surjections $q_1:\hat{\mathcal E}\to \Graph_1$ and $q_2:\hat{\mathcal E}\to \Graph_2$, the set-valued map 
\[
Q_{12}:q_1\circ q_2^{-1}:\Graph_2\rightrightarrows \Graph_1
\]
 allows us to identify points in the metric graphs $\Graph_1:=\faktor{\hat{\mathcal E}}{\approx_1}$ and $\Graph_2:=\faktor{\hat{\mathcal E}}{\approx_2}$. 
% In the following, we always let 
%$\Graph:=\faktor{\mathcal E}{ \sim}=\bigsqcup\limits_{\me\in\mE} [0,\ell_\me]/ \sim$ be a (finite or countably infinite) metric graph.

Given two metric graphs $\Graph_1,\Graph_2$, we write $\Graph_1\equiv \Graph_2$ if both $\Graph_1,\Graph_2$ are subdivisions of the same metric graph $\Graph$. Now, $\equiv$ is an equivalence relation on the set of all metric graphs: $\Graph_1\equiv \Graph_2$ if they agree modulo removing vertices of degree 2.

\begin{defi}
Let $\Graph$ be a metric graph.
We call \textit{primitive metric graph} of $\Graph$ the equivalence class $\Graphprim=[\Graph]$ with respect to the above defined equivalence relation $\equiv$.
\end{defi}

\begin{defi}
Let $\Graph$ be a metric graph.
\begin{enumerate}[(i)]
\item Given $\mE_0\subset \mE$, consider $\mathcal E_0:=\bigsqcup\limits_{\me\in \mE_0}[0,\ell_\me]$ and the set $\mathcal V_0:=\bigsqcup\limits_{\me\in\mE_0}\{0,\ell_\me\}$ of endpoints of $\mathcal E_0$. The metric graph $\Graph_0:=\faktor{ \mathcal E_0}{\sim_0}$ is said to be a \textit{metric subgraph} of $\Graph:=\faktor{\mathcal E}{\sim}$, where $\sim_0$ is the restriction of $\sim$ to $\mathcal V_0$. 
\item With a slight abuse, we also say that $\Graph_0$ is a metric subgraph of $\Graph$ if any representative of $[\Graph_0]$ is a metric subgraph, in the sense of (i), of any representative of $[\Graph]$.

\item A \textit{connected component} of a metric graph $\Graph$ is a metric subgraph $\Graph_0$ of $\Graph$ that is maximal (wrt to $\subset$ for $\mE_0$) among connected ones.
\end{enumerate}
\end{defi}

\section{Function spaces}\label{sec:functsp}

 Whenever considering a \textit{continuous} function $f$ on a metric graph $\Graph$, there is a uniquely determined continuous function induced by $f$ on any further metric graph belonging to $\Graphprim=[\Graph]$. It would be appropriate to consider the space of continuous functions $C(\Graphprim)$, yet in practice the notation $C(\Graph)$ is customary in the literature: this space is isometrically isomorphic to the space of continuous functions supported on any other representative of $[\Graphprim]$.

Similarly, two functions on $\Graph$ can be identified if they agree up to a Lebesgue null set. Accordingly, any measurable $f:\Graph\to \K$ can -- up to a Lebesgue null set -- be canonically identified with a unique function defined on any rearrangement of $\Graph$: accordingly, the Lebesgue space $L^p(\Graph)$ is isomorphic to  $L^p(\Graph')$ for any $p\in [1,\infty]$.

Summing up, we can introduce the function spaces 
\[
C(\Graph)\quad\hbox{and}\quad L^p(\Graph),\qquad 1\le p\le\infty
\]
and then, recursively, for all $k\in\mathbb N$ the spaces
\[
C^k(\Graph):=\left\{f=\bigoplus\limits_{\me\in\mE}f_\me \in \bigoplus\limits_{\me\in\mE}C^{k}([0,\ell_\me]):f^{(h)}:=\bigoplus\limits_{\me\in\mE}f^{(h)}_\me\in C(\Graph)\hbox{ for all }1\le h\le k\right\}
\]
and
\[
\begin{split}
W^{k,p}(\Graph)&:=\{f\in L^p(\Graph):f^{(h)}\in C(\Graph)\hbox{ for all }0\le h\le k-1\\
&\qquad \qquad\qquad \hbox{ and }f^{(j)}\in L^p(\Graph)\hbox{ for all }0\le j\le k\},\qquad 1\le p\le\infty.
\end{split}
\]
Observe that $W^{k,p}(\Graph)$ is the closure of $C^{k}(\Graph)$ with respect to the norm $\|f\|_{k,p}:=\sum\limits_{h=0}^k \|f^{(h)}\|_p$. 

%If $\mathcal G$ has finite volume, then Observe that $W^{k,p}(\Graph)$ is canonically identified with a subspace of $\bigoplus_{\me\in\mE} W^{k,p}(0,\ell_\me)$, which is canonically identified with $L

\section{Graph operations} 
If two metric graphs $\Graph_1,\Graph_2$ are defined upon the same $\mathcal E$, they are completely characterized by the equivalence relations $\sim_1,\sim_2$. Accordingly, we can easily define binary operations on metric graphs by means of operations involving $\sim_1,\sim_2$. Recalling that given any binary relation $A\subset \mathcal V\times \mathcal V$, the \textit{equivalence relation generated by $A$} is by definition the intersection of the equivalence relations on $\mathcal V$ that contain $A$, we can, e.g., consider
\begin{itemize}
\item the \textit{intersection} of  $\Graph_1,\Graph_2$ is the metric graph on $\mathcal E$ obtained by taking $\sim$ to be $\sim_1\cap \sim_2$ (this is automatically an equivalence relation!); 
\item the \textit{union} of  $\Graph_1,\Graph_2$ is the metric graph on $\mathcal E$ obtained by taking $\sim$ to be the equivalence relation generated by $\sim_1\cup \sim_2$ (the latter is automatically reflexive and symmetric).
\end{itemize}

%Unfortunately, these notions aren't stable under insertion of dummy vertices.
% While the definitions of union and intersection can probably be patched, the definition of complement certainly cannot.

\begin{example}
Take $\Graph_1$ to be a cycle consisting of two edges; and $\Graph_2$ to be the disconnected graph consisting of two loops, each consisting of one edge. In the above formalism, they are modeled by taking $\mE=\{1,2\}$ and, for any $\ell_1,\ell_2\in (0,\infty)$, by the equivalence relations
\[
\begin{tabular}{c|c|c|c|c|}
& $(0,1)$ & $(\ell_1,1)$ & $(0,2)$ & $(\ell_2,2)$ \\ 
\hline 
$(0,1)$ & $\times$ & &  & $\times$\\ 
\hline 
$(\ell_1,1)$ &  & $\times$ & $\times$ &  \\ 
\hline 
 $(0,2)$ &  & $\times$ & $\times$ &  \\ 
\hline 
$(\ell_2,2)$ &  $\times$ &  & & $\times$ \\ 
\hline 
\end{tabular} 
\qquad \hbox{and}\qquad 
\begin{tabular}{c|c|c|c|c|}
& $(0,1)$ & $(\ell_1,1)$ & $(0,2)$ & $(\ell_2,2)$ \\ 
\hline 
$(0,1)$ & $\times$ & $\times$ &  & \\ 
\hline 
$(\ell_1,1)$ & $\times$ & $\times$ &  &  \\ 
\hline 
 $(0,2)$ &  &  & $\times$ & $\times$ \\ 
\hline 
$(\ell_2,2)$ &  &  &$\times$ & $\times$ \\ 
\hline 
\end{tabular} 
\]
respectively. Their intersection and union are given by the equivalence relations
\[
\begin{tabular}{c|c|c|c|c|}
& $(0,1)$ & $(\ell_1,1)$ & $(0,2)$ & $(\ell_2,2)$ \\ 
\hline 
$(0,1)$ & $\times$ & &  & \\ 
\hline 
$(\ell_1,1)$ &  & $\times$ &  &  \\ 
\hline 
 $(0,2)$ &  & & $\times$ &  \\ 
\hline 
$(\ell_2,2)$ & &  & & $\times$ \\ 
\hline 
\end{tabular} 
\qquad \hbox{and}\qquad 
\begin{tabular}{c|c|c|c|c|}
& $(0,1)$ & $(\ell_1,1)$ & $(0,2)$ & $(\ell_2,2)$ \\ 
\hline 
$(0,1)$ & $\times$ & $\times$ & $\times$ & $\times$ \\ 
\hline 
$(\ell_1,1)$ & $\times$ & $\times$ & $\times$ & $\times$ \\ 
\hline 
 $(0,2)$ & $\times$ & $\times$ & $\times$ & $\times$ \\ 
\hline 
$(\ell_2,2)$ & $\times$ & $\times$ & $\times$ & $\times$ \\ 
\hline 
\end{tabular} 
\]
respectively, i.e., they correspond to two disjoint intervals and to the figure-8 graph, respectively.
\end{example}

\begin{remark}
Following the above path, we can also define the \textit{complement} of  $\Graph_2$ in $\Graph_1$ as the metric graph on $\mathcal E$ obtained by taking $\sim$ to be the equivalence relation generated by $\sim_1\setminus \sim_2$; and, canonically, the \textit{complement of $\Graph:=\Graph_2$} obtained by taking $\Graph_1$ to be the flower graph (much like in the combinatorial graph setting, where the canonical ambient graph is the complete one).

Complements of metric graph tend to be trivial, though. Take e.g.\ a lasso graph: formally, it is given by $\mE=\{1,2\}$ and, for any $\ell_1,\ell_2\in (0,\infty)$, by the equivalence relation
\[
\begin{tabular}{c|c|c|c|c|}
& $(0,1)$ & $(\ell_1,1)$ & $(0,2)$ & $(\ell_2,2)$ \\ 
\hline 
$(0,1)$ & $\times$ & $\times$ & $\times$ & \\ 
\hline 
$(\ell_1,1)$ & $\times$ & $\times$ & $\times$ &  \\ 
\hline 
 $(0,2)$ & $\times$ & $\times$ & $\times$ &  \\ 
\hline 
$(\ell_2,2)$ &  &  & & $\times$ \\ 
\hline 
\end{tabular} 
\]
 on $\{(0,1),(\ell_1,1),(0,2),(0,\ell_2)\}$; the equivalence relation generated by its complement yields
\[
\begin{tabular}{c|c|c|c|c|}
& $(0,1)$ & $(\ell_1,1)$ & $(0,2)$ & $(\ell_2,2)$ \\ 
\hline 
$(0,1)$ & $\times$ & $\times$ & $\times$ & $\times$ \\ 
\hline 
$(\ell_1,1)$ & $\times$ & $\times$ & $\times$ & $\times$ \\ 
\hline 
 $(0,2)$ & $\times$ & $\times$ & $\times$ & $\times$ \\ 
\hline 
$(\ell_2,2)$ & $\times$ & $\times$ & $\times$ & $\times$ \\ 
\hline 
\end{tabular} 
\]
i.e., the complement of the lasso graph is the figure-8 graph. Likewise,  the figure-8 graph is also the complement of the cycle (formally consisting of two edges) as well as complement of the disconnected graphs consisting of either two intervals or of two loops.
\end{remark}

\bibliographystyle{alpha}
\bibliography{../../referenzen/literatur}

\end{document}